%

\magnification=\magstep1
\input amstex
\UseAMSsymbols
\input pictex
\vsize=23truecm
\NoBlackBoxes
\parindent=20pt
   
   \font\rmk=cmr8    \font\itk=cmti8  \font\ttk=cmtt8


\def\NC{\operatorname{NC}}
\def\mod{\operatorname{mod}}

\def\Hom{\operatorname{Hom}}
\def\End{\operatorname{End}}

\def\op{{\text{op}}}
\def\Ext{\operatorname{Ext}}

\def\NC{\operatorname{NC}}
\def\bdim{\operatorname{\bold{dim}}}

\def\Rahmen#1%
   {$$\vbox{\hrule\hbox%
                  {\vrule%
                       \hskip0.5cm%
                            \vbox{\vskip0.3cm\relax%
                               \hbox{$\displaystyle{#1}$}%
                                  \vskip0.3cm}%
                       \hskip0.5cm%
                  \vrule}%
           \hrule}$$}


\vglue1.5truecm
\centerline{\bf The number of complete exceptional sequences for a Dynkin algebra}
	\bigskip\medskip
\centerline{Mustafa A\. A\. Obaid, S\. Khalid Nauman, Wafa S\. Al Shammakh,}
	\smallskip
\centerline{Wafaa M\. Fakieh and Claus Michael Ringel (Jeddah)}
	\bigskip\bigskip

\plainfootnote{}
{\rmk 2010 \itk Mathematics Subject Classification. \rmk 
Primary:
        16G20,  
        16G60, 
        05A19, 
        05E10. 
Secondary:
        16D90, 
        16G70. 
        16G10. 
     Lattice of non-crossing partitions. Binomial convolution. Abel's identity.
     Categorification. 
}

{\narrower\narrower Abstract. The Dynkin algebras are the 
hereditary artin algebras of finite representation type. The paper determines the 
number of complete exceptional sequences for any Dynkin algebra. 
Since the complete exceptional sequences for a Dynkin algebra of Dynkin type $\Delta$
correspond bijectively to the maximal chains in the lattice of non-crossing partitions 
of type $\Delta$, the calculations presented here may also be considered as a
categorification of the corresponding result for non-crossing partitions.
\par}
	\bigskip\bigskip 
{\bf 1\. Introduction.}
	\medskip 
We consider Dynkin algebras $\Lambda$, 
these are the hereditary artin algebras of finite representation type. 
Note that the indecomposable $\Lambda$-modules correspond
bijectively to the positive roots of a Dynkin diagram 
$\Delta(\Lambda)$; such a diagram is the disjoint union of connected diagrams
and the connected Dynkin diagrams are of the form $\Bbb A_n, \Bbb B_n,\dots, \Bbb G_2$.
Let us remark that the vertices $i$
of $\Delta(\Lambda)$ correspond bijectively to the simple 
$\Lambda$-modules, there is an edge between two vertices if and only if there is a
non-trivial extension between the corresponding simple modules (in one of the two possible directions), and the
lacing (in the cases $\Bbb B_n, \Bbb C_n, \Bbb F_4, \Bbb G_2$) records the relative size of the endomorphism rings of the simple modules, see [DR1] or [DR2].
We call $\Lambda$ a {\it Dynkin algebra} of type $\Delta(\Lambda),$ the number of simple
$\Lambda$-modules will be called the {\it rank} of $\Lambda$
(let us stress the following: when we refer to the number of modules of some kind or the
number of sequences of modules, then we mean of course the number of isomorphism classes).

Given a Dynkin algebra $\Lambda$
an {\it exceptional sequence} for $\Lambda$  is a sequence
$(M_1,\dots, M_t)$ of indecomposable $\Lambda$-modules 
such that $\Hom(M_i,M_j) = 0 = \Ext^1(M_i,M_j)$ for $i > j.$ 
 The cardinality of an exceptional sequence
is bounded by the rank $n$ of $\Lambda$ and 
the exceptional sequences of cardinality $n$ are said to be {\it complete}.
Any exceptional sequence $(M_1,\dots,M_t)$
can be extended to a complete exceptional sequence $(M_1,\dots,M_n)$; 
in case $t = n-1$, the extension is unique (for all these assertions, see [CB] and [R2]). 

Let $e(\Lambda)$ be the number 
of complete exceptional sequences for the Dynkin algebra 
$\Lambda$. 
In case $\Lambda$ is the path algebra of a quiver, the number $e(\Lambda)$ has
been determined by Seidel [Se] in 2001. The aim of this note is to finalize these investigations by
dealing also with the Dynkin diagrams which are not simply laced.
There are direct connections between the representation theory of a Dynkin algebra 
$\Lambda$  and the lattice $L$ of non-crossing partitions
of type $\Delta(\Lambda)$ which we will outline at the end of the introduction. 
In particular, the complete exceptional sequences for $\Lambda$ correspond bijectively to the
maximal chains in $L.$ 
Thus, the calculations may also be considered as a categorification 
of the corresponding result for $L$. 

As we will see, the number $e(\Lambda)$ 
only depends on $\Delta = \Delta(\Lambda)$, thus we may write
$e(\Delta)$ instead of $e(\Lambda)$. Also, the shuffle lemma 
presented in section 2 shows that
it is sufficient to look at the connected Dynkin diagrams $\Delta$. 
	\bigskip 
The following table exhibits the numbers $e(\Delta)$ for 
any connected Dynkin diagram $\Delta$:
$$
{\beginpicture
\setcoordinatesystem units <1.5cm,.7cm>
\plot -.5 0.5  8.2 0.5 /
\plot 0.3 1.2 0.3 -.3 /
\put{$\Delta$} at -.2 1
\put{$e(\Delta)$} at -.2 0
\put{$\Bbb A_n$} at 1 1
\put{$\Bbb B_n,\Bbb C_n$} at 2 1
\put{$\Bbb D_n$} at 3 1
\put{$\Bbb E_6$} at 4 1
\put{$\Bbb E_7$} at 5 1
\put{$\Bbb E_8$} at 6 1
\put{$\Bbb F_4$} at 7 1
\put{$\Bbb G_2$} at 7.9 1
\put{$(n\!+\!1)^{n-1}$\strut} at 1 0
\put{$n^{n}$\strut} at 2 0
\put{$2(n\!-\!1)^{n}$\strut} at 3 0
\put{$2^9\!\cdot\! 3^4$\strut} at 4 0
\put{$2\!\cdot\! 3^{12}$\strut} at 5 0
\put{$2\!\cdot\! 3^5\!\cdot\! 5^7$\strut} at 6 0
\put{$2^4\!\cdot\! 3^3$\strut} at 7 0
\put{$2\!\cdot\! 3$\strut} at 7.9 0

\endpicture}
$$
It seems to be of interest that the numbers $e(\Delta)$ have only few different prime factors, 
all of them being rather small.
Using the table, one easily verifies the following
remarkable formula
\Rahmen{e(\Delta) = \frac{n!\,h(\Delta)^n}{|W(\Delta)|}}
where $W(\Delta)$ is the Weyl group of type $\Delta$ and $h(\Delta)$ the
corresponding Coxeter number. Here are the numbers in question, as given, for example, 
in the appendix of [B]:
$$
{\beginpicture
\setcoordinatesystem units <1.5cm,.7cm>
\plot -.5 0.5  8.2 0.5 /
\plot 0.3 1.2 0.3 -1.3 /
\put{$\Delta$} at -.2 1
\put{$h(\Delta)$} at -.2 0
\put{$|W(\Delta)|$} at -.2 -1
\put{$\Bbb A_n$} at 1 1
\put{$\Bbb B_n,\Bbb C_n$} at 2 1
\put{$\Bbb D_n$} at 3 1
\put{$\Bbb E_6$} at 4 1
\put{$\Bbb E_7$} at 5 1
\put{$\Bbb E_8$} at 6 1
\put{$\Bbb F_4$} at 7 1
\put{$\Bbb G_2$} at 7.9 1
\put{$n\!+\!1$\strut} at 1 0
\put{$2n$\strut} at 2 0
\put{$2(n\!-\!1)$\strut} at 3 0
\put{$2^2\!\cdot\! 3$\strut} at 4 0
\put{$2\!\cdot\! 3^{2}$\strut} at 5 0
\put{$2\!\cdot\! 3\!\cdot\! 5$\strut} at 6.1 0
\put{$2^2\!\cdot\! 3$\strut} at 7 0
\put{$2\!\cdot\! 3$\strut} at 7.9 0

\put{$(n\!+\!1)!$\strut} at 1 -1
\put{$2^nn!$\strut} at 2 -1
\put{$2^{n-1}n!$\strut} at 3 -1
\put{$2^73^45$\strut} at 4 -1
\put{$2^{10}3^{4}5\!\cdot\! 7$\strut} at 5 -1
\put{$2^{14}3^55^27$\strut} at 6.1 -1
\put{$2^7\!\cdot\! 3^2$\strut} at 7.05 -1
\put{$2^2\!\cdot\! 3$\strut} at 7.9 -1

\endpicture}
$$
Unfortunately, our proof does
not provide any illumination of the formula 
(and we should admit that the observation that the formula holds is 
stolen from Chapoton [Ch], see the end of the introduction).

As we have mentioned, for $\Lambda$ the path algebra of a quiver (thus for the typical Dynkin
algebras of type $\Bbb A_n, \Bbb D_n, \Bbb E_n$), the numbers $e(\Lambda)$
have been determined already by Seidel [Se] in 2001. The essential cases which
were missing are the Dynkin algebras of type $\Bbb B_n$. The inductive strategy of proof
works for all types. However, we also will show a direct relationship between the cases
$\Bbb B_n$ and $\Bbb A_{n-1}$, and this could be used directly in order to complete Seidel's
considerations. Clearly, for $n=2$, the number $e(\Lambda)$ is just the number of indecomposable
modules, in particular we have $e(\Bbb G_2) = 6.$
	\medskip 
Here is an outline of the proof: we will use induction on the rank $n$ of $\Lambda$.
If $M$ is an indecomposable $\Lambda$-module, 
$M^\perp$ be the full subcategory of $\mod\Lambda$ consisting of
all modules $N$ such that $\Hom(M,N) = 0 = \Ext^1(M,N).$ 
Since $M$ is exceptional, one knows that $M^\perp$ 
is  (equivalent to) the 
module category of a hereditary artin algebra of rank $n-1$ (see [GL] or [S2]),
thus by induction we may assume to know $e(M^\perp).$ 
Obviously, the 
complete exceptional sequences $(M_1,\dots,M_n)$ with $M_n = M$ correspond bijectively to the 
complete exceptional sequences in  $M^\perp$, thus $e(M^\perp)$ is the number of complete exceptional
sequences in $\mod\Lambda$ whose last entry is $M$. In section 3 we will see that there is a 
vertex $i_M$ of $\Delta$ such that 
$e(M^\perp) = e(\Delta(i_M))$, where $\Delta(i)$ is obtained from $\Delta = 
\Delta(\Lambda)$ by deleting the vertex $i$
and all the edges involving $i$.
Thus
$$
 e(\Delta) = \sum\nolimits_{M} e(\Delta(i_M)),
$$
and therefore, for $\Delta$ being connected, there is the following reduction formula
$$
 e(\Delta) = \frac h2 \sum\nolimits_{i\in \Delta_0} e(\Delta(i))
$$
where $h$ is the Coxeter number for $\Delta$ (see section 4). 
In section 5 we will use the reduction formula in order to obtain 
the entries of the table, here we have to proceed case by case.
The proof of cases $\Bbb A_n, \Bbb B_n, \Bbb C_n, \Bbb D_n$ relies on some 
well-known recursion formulas which go back to Abel [Ab], see the Appendix. 
Conversely, one may observe that
the interpretation using complete exceptional sequences provides a categorification
of these formulas. 

Since we deal with artin algebras (and not more generally with artinian
rings), the diagrams which arise are the Dynkin diagrams $\Bbb A_n,\dots \Bbb G_2$. 
  If one is interested in all the finite Coxeter diagrams (thus also in $\Bbb I_2(m), \Bbb H_3, \Bbb H_4$),
  one may consider in the same way corresponding artinian rings
(they are known to exist
for $\Bbb I_2(5), \Bbb H_3, \Bbb H_4$, see [S1] as well as [DRS] and [O]), 
this will be done in [FR].

	\bigskip
{\bf The general frame.}
The calculations presented here can be seen in a broader frame,
since the representation theory
of hereditary artinian rings has turned out to be an intriguing tool for dealing with
various questions in different parts of mathematics.
In particular, there is a strong relationship to the theory
of (generalized) non-crossing partitions (see for example [Ag]) as observed first by Fomin and Zelevinsky.
As Ingalls and Thomas [IT] have shown, given the path algebra $\Lambda$
of a finite directed quiver of type $\Delta$,
there is a poset isomorphism between the poset of thick subcategories of $\mod\Lambda$ with 
generators and the poset $\NC(\Delta)$ of non-crossing partitions of type $\Delta$ (and this result can
easily  be extended to arbitrary hereditary artin algebras $\Lambda$); we recall that
a full subcategory is said to be {\it thick} (or ``wide'') provided it is closed 
under kernels, cokernels and extensions. Of course, in case $\Lambda$ is of finite
representation type, any thick subcategory has a generator. Hubery and Krause [HK] have
pointed out that the Ingalls-Thomas bijection yields a bijection between the complete exceptional 
sequences for $\Lambda$ and the maximal chains in the poset $\NC(\Delta)$.
Namely, given a complete exceptional sequence $(M_1,\dots,M_n)$ for $\Lambda$
let $\Cal U_i = (M_{i+1}\oplus \cdots
\oplus M_n)^\perp$, for $0\le i \le n$. 
Then $0 = \Cal U_0 \subset \Cal U_1 \subset \cdots \subset \Cal U_n = \mod\Lambda$ is a maximal chain
of thick subcategories of $\mod\Lambda$ with generators. Conversely, let us assume that 
$0 = \Cal U_0 \subset \Cal U_1 \subset \cdots \subset \Cal U_n = \mod\Lambda$ is a maximal chain of thick subcategories
of $\mod\Lambda$ with generators. Then $\Cal U_{n-1}$  is the module category of a hereditary artin algebra of
rank $n-1$, thus 
by induction the chain $0 = \Cal U_0 \subset \Cal U_1 \subset \cdots \subset \Cal U_{n-1}$ corresponds to a complete exceptional sequence $(M_1,\dots,M_{n-1})$ in $\Cal U_{n-1}$, and this is an exceptional sequence for $\Lambda$ 
of cardinality $n-1$. As we have mentioned, 
there is a uniquely determined $\Lambda$-module $M_n$ such that
$(M_1,\dots,M_n)$ is a complete exceptional sequence for $\Lambda$. We see that there
is a canonical bijection between the complete exceptional sequences for $\Lambda$ and the
set of maximal chains of thick subcategories of $\mod\Lambda$ with generators,
thus with the maximal chains in $\NC(\Delta)$. 

This shows that the numbers $e(\Delta)$ calculated here for the Dynkin diagrams $\Delta$
via representation theory are nothing else
than the numbers of maximal chains in $\NC(\Delta)$ (in the Dynkin case,
this poset is even a lattice) or, equivalently, the numbers of  
factorizations of a fixed Coxeter element 
as a product of $n$ reflections. The latter numbers for 
$\Delta = \Bbb A_n, \Bbb B_n, \Bbb D_n$
have been determined in a famous letter [D] of Deligne to Looijenga.
The numbers of maximal chains in $\NC(\Delta)$ have been
calculated for the cases $\Bbb A_n,$ $\Bbb B_n$ and $\Bbb D_n$ by Kreweras [K], 
Reiner [Rn] and Athanasiadis-Reiner [AR], respectively, and in general by 
Chapoton [Ch] and Reading [Rd], see also Chapuy-Stump [CS].
It seems that the term $n!h^n/|W|$ is mentioned first by Chapoton [Ch].

The present paper only relies on well-known properties of the module category of an artin
algebra. On the other hand, the result presented here, and indeed also the main steps 
of our proof, may be considered as a categorification
of the considerations of Deligne and Reading. 

The authors are strongly indebted to 
H\. Krause, C\. Stump and H\. Thomas for pointing out pertinent
references concerning non-crossing partitions and the relevance of the numbers 
$e(\Delta)$, and to M\. Baake for helful remarks concerning the binomial convolution.
The references [AR], [Rn] were provided by Thomas, the references [D], [CS] and [Rd]
by Krause. Also, 
we learned from Krause that in the context of  simple singularities,
the numbers $e(\Delta)$ for simply laced Dynkin diagrams $\Delta$
have been presented in 1974 by Looijenga [L].

	\bigskip
{\bf Acknowledgment.}
This work is funded by the Deanship of Scientific Research, 
King Abdulaziz University, under grant No. 2-130/1434/HiCi. 
The authors, therefore, acknowledge technical and financial support of KAU.
	\bigskip\medskip
{\bf 2. The shuffle lemma.}
	\medskip
{\bf Lemma 1 (Shuffle Lemma).} {\it Let $\Lambda, \Lambda'$ be representation-finite hereditary artin algebras of ranks 
$n,n'$
respectively. Then}
$$
 e(\Lambda\times\Lambda') = \binom{n+n'}n
 e(\Lambda) e(\Lambda').
$$

Proof. Let $(E_1,\dots,E_n)$ be a (complete) exceptional sequence in $\mod\Lambda$ and let
$(E'_1,\dots,E'_{n'})$ be a (complete) exceptional sequence in $\mod\Lambda'$. Let $I$ be a subset of 
$\{1,2,\dots,n+n'\}$ of cardinality $n$, say let $I = \{i_1< i_2 < \dots < i_n\}$  and let
$\{j_1 < j_2 < \cdots < j_{n'}\}$ be its complement. 
Let $(M_1,\dots, M_{n+n'})$ be defined by $M_{i_t} = E_t$ for $1\le t \le n$ and 
$M_{j_t} = E'_t$ for $1\le t \le n'.$ Then clearly $(M_1,\dots, M_{n+n'})$ is a complete exceptional sequence
in $\mod(\Lambda\times\Lambda')$ and every  complete exceptional sequence
in $\mod(\Lambda\times\Lambda')$ is obtained in this way. Thus, fixing a subset $I$ of cardinality $n$, 
the number of  complete exceptional sequences
$(M_1,\dots, M_{n+n'})$  in $\mod(\Lambda\times\Lambda')$ with $M_i$ in $\mod\Lambda$ for all $i\in I$
is equal to $e(\Lambda) e(\Lambda')$, and the number of such subsets $I$ is just $\binom{n+n'}n$.
This completes the proof. 
	\bigskip\bigskip

{\bf 3. The category $M^\perp$.}
	\medskip 
Let $\Lambda$ be a representation-finite hereditary artin algebra of rank $n$. Let $\Delta =
\Delta(\Lambda)$. Given a vertex $i$ of $\Delta,$ let $\Delta(i)$ be obtained from $\Delta$ by deleting
the vertex $i$ and the edges involving $i$ (it is of course again a Dynkin diagram).
	\medskip
Let $\tau$ be the Auslander-Reiten translation for $\Lambda$. For every indecomposable $\Lambda$-module $M$,
there is a natural number $t$ such that $\tau^tM$ is indecomposable projective, thus $\tau^tM = P(i_M)$ for a
(uniquely determined) vertex $i_M$ of $\Delta.$ 
	\medskip 
Let $M$ be an indecomposable module. It is known that the category $M^\perp$ is equivalent to
a module category $\mod\Lambda'$ where $\Lambda'$ is a representation-finite hereditary artin algebra of rank $n-1$.
	\medskip 
{\bf Lemma 2.} {\it Let $M$ be an indecomposable module and assume that $M^\perp$ is equivalent to
the module category $\mod\Lambda'$. Then $\Lambda'$ has type $\Delta(i_M)$.} 
	\medskip 
Proof. First, assume that $M = P(i)$ is indecomposable projective, thus $i = i_M$.
Let $\epsilon_i$ be an idempotent of $\Lambda$ such that $P(i) = \Lambda \epsilon_i.$ Then $M^\perp$
is the set of $\Lambda$-modules $N$ with $\Hom(P(i),N) = 0,$ thus the set of $\Lambda/\Lambda \epsilon_i\Lambda$-modules.
On the other hand, we have $\Delta(\Lambda/\Lambda \epsilon_i\Lambda) = \Delta(i)$. 

Now assume that $M$ is indecomposable and not projective. 
There is a slice $\Cal S$ (in the sense of [R2]) in the Auslander-Reiten quiver
of $\Lambda$ such that $M$ is a sink for $\Cal S$. Let $M_1,\dots M_n$ be the indecomposable modules in $\Cal S$,
one from each isomorphism class, and we assume that $M_n = M.$ 
Since $M$ is a sink of $\Cal S$, we know that $\Hom(M,M_i) = 0$ for $1\le i \le n-1$, thus
the modules $M_1,\dots, M_{n-1}$ belong to $M^\perp.$ 
Let $T = \bigoplus_{i=1}^{n-1} M_i$, then $T$ is a tilting module for $M^\perp = \mod\Lambda'$ 
(it has no self-extensions
and enough indecomposable direct summands). Since $\Cal S$ is a slice, we know that the endomorphism ring of
$\bigoplus_{i=1}^n M_i$ is hereditary, thus also $\End(T)^\op$ is 
hereditary and the Dynkin diagram $\Delta(\End(T)^\op)$ is just $\Delta(i_M)$. 
A tilting module with hereditary endomorphism ring is a slice module
(see for example [R3], section 1.2). 
Thus $T$ is a slice module for $\mod\Lambda'$ and therefore $\Lambda'$ 
and $\End(T)^\op$ have the same Dynkin type.  
This shows that the Dynkin type of $\Lambda'$ is $\Delta(i_M)$.
	\bigskip\bigskip 

{\bf 4\. The reduction formula.}
	\medskip 
We assume by induction that $e(\Lambda')$ only depends on $\Delta(\Lambda')$ for any 
representation-finite hereditary artin algebra $\Lambda'$ of rank $n' < n$. 
	\medskip 
{\bf Proposition.} {\it Let $\Lambda$ be a connected 
representation-finite hereditary artin algebra of rank $n$
and type $\Delta$. Then
$$
 e(\Lambda) = \frac h2 \sum\nolimits_{i\in \Delta_0} e(\Delta(i)),
$$
where $h$ is the Coxeter number for $\Delta$.}
	\medskip
\noindent
This reduction formula shows that $e(\Lambda)$ only depends on $\Delta = \Delta(\Lambda).$
	\bigskip 
Proof. If $M$ is an indecomposable $\Lambda$-module, then we have seen in section 3 that
$M^\perp$ is equivalent to the module category $\mod\Lambda'$, where $\Lambda'$ is of type $\Delta(i_M)$.
Thus  
$$
 e(M^\perp) = e(\Delta(i_M)).
$$

For any vertex $i$ of $\Delta,$ let $m(i)$ be the length of the $\tau$-orbit of $P(i)$, thus there
are precisely $m(i)$ indecomposable modules $M$ such that $i_M = i.$ Therefore
$$
 e(\Lambda) = \sum\nolimits_{M} e(M^\perp) =\sum\nolimits_{M} e(\Delta(i_M)) = \sum\nolimits_i m(i)e(\Delta(i)).
$$

We have to distinguish two cases. First, assume that $\Delta$ is not of the form $\Bbb A_n$ or 
$\Bbb D_{2m+1}$ or $\Bbb E_6$. In this case, we have $m(i) = \frac h2$ for any vertex $i$ of $\Delta$. Therefore
$$
 \sum\nolimits_i m(i)e(\Delta(i)) = 
 \sum\nolimits_i \frac h2 e(\Delta(i)).
$$

Second, assume that $\Delta$ is equal to $\Bbb A_n$, or 
$\Bbb D_{2m+1}$ or $\Bbb E_6$. Thus, there is a (unique) automorphism $\rho$ of $\Delta$
of order 2. One knows that $m(i)+m(\rho(i)) = h$ for all vertices $i$ of $\Delta$. The automorphism $\rho$ shows that
$e(\Delta(\rho(i))) = e(\Delta(i))$, thus 
$$
\align
 2\sum\nolimits_i m(i)e(\Delta(i)) &= \sum\nolimits_i m(i)e(\Delta(i)) + \sum\nolimits_i m(\rho(i))e(\Delta(\rho(i))) \cr
 &= \sum\nolimits_i (m(i)+m(\rho(i))e(\Delta(i)) \cr
 &=\sum\nolimits_i h \cdot e(\Delta(i)).
\endalign
$$
Dividing by $2$ we obtain the required formula.

	\bigskip\bigskip
{\bf 5. The different cases.}
	\medskip
{\bf Type $\Bbb A_n$.} This concerns the following diagram
$$
{\beginpicture
\setcoordinatesystem units <1cm,1cm>
\multiput{$\circ$} at 0 0  1 0  2 0  4 0 /
\put{$\cdots$} at 3 0
\plot 0.1 0  0.9 0 /
\plot 1.1 0  1.9 0 /
\plot 2.1 0  2.6 0 /
\plot 3.4 0  3.9 0 /
\put{$0$} at 0 -.3
\put{$1$} at 1 -.3
\put{$2$} at 2 -.3
\put{$n\!-\!1$} at 4 -.3
\endpicture}
$$
We have $\Delta(i) = \Bbb A_{i}\sqcup \Bbb A_{n-i-1}$, therefore, by the shuffle lemma and induction,
$$
 e(\Delta(i)) = \tbinom{n-1}{i} e(A_{i}) e(\Bbb A_{n-i-1}) = \tbinom{n-1}{i} (i+1)^{i-1}(n-i)^{n-i-2}.
$$
Thus we have to calculate
$$
 \sum\nolimits_{i=0}^{n-1} e(\Delta(i)) = 
 \sum\nolimits_{i=0}^{n-1}\tbinom{n-1}{i} (i+1)^{i-1}(n-i)^{n-i-2}, 
$$
but this is the coefficient $F(n-1)$ of the power series $F = A*A$, see the appendix,
and the formula (1) asserts that $F(n\!-\!1) = 2(n+1)^{n-2}.$

Now $h = n+1,$ thus
$$
 \frac h2 \sum\nolimits_{i=1}^n e(\Delta(i)) = \frac {n+1}2 2 (n+1)^{n-2} = (n+1)^{n-1}.
$$

	\bigskip\bigskip
\vfill\eject
{\bf Type $\Bbb B_n:$ The relationship between $\Bbb B_n$ and $\Bbb A_{n-1}$.}
	\medskip
Let us directly show the following relationship:
$$
 e(\Bbb B_n) = n^2\cdot e(\Bbb A_{n-1}).
$$
	\medskip
Proof. Let $\Lambda$ be a hereditary artin algebra of type $\Bbb B_n$. Let $P$ be the indecomposable 
projective $\Lambda$-module such that
$\bdim P$ is a short root. 
If $(M_1,\dots, M_n)$ is an exceptional sequence in $\mod\Lambda$, then there is precisely one index
$i$ such that $\bdim M_i$ is a short root (see [R2]). Thus, let $\Cal E_i(\mod\Lambda)$ be the set 
of exceptional sequences in $\mod\Lambda$ such that $\bdim M_i$ is a short root, and let $e_i(\mod\Lambda)$
the cardinality of $\Cal E_i(\mod\Lambda)$. 
If $i < n$, and $(M_1,\dots, M_n)$ belongs to $\Cal E_i(\mod\Lambda)$, then there is a uniquely
determined  element $(M_1,\dots, M_{i-1},M_{i+1},M_i^*,M_{i+2},\dots, M_n)$ in $\Cal E_{i+1}(\mod\Lambda)$
and every element of $\Cal E_{i+1}(\mod\Lambda)$ is obtained in this way (again, see [R2]). 
This shows that
$e_i(\mod\Lambda) = e_{i+1}(\mod\Lambda)$ 
and  therefore 
$$
 e(\Lambda) = \sum\nolimits_{i=1}^n e_i(\Lambda) = n\cdot e_n(\Lambda).
$$
There are precisely $n$ indecomposable modules $M$ such that $\bdim M$ is a short root, namely the
modules in the $\tau$-orbit $\Cal O(P)$ of $P$. For any module $M$ in $\Cal O(P)$, the exceptional
sequences $(M_1,\dots,M_n)$ with $M_n = M$ correspond bijectively to the exceptional sequences
in $M^\perp$, and $M^\perp$ is equivalent to a module category $\mod\Lambda_M$ with $\Lambda_M$ a hereditary
artin algebra of type $\Bbb A_{n-1}$. This shows that
$$
 e_n(\mod\Lambda) = \sum\nolimits_{M\in \Cal O(P)} e(M^\perp) = n\cdot e(\Bbb A_{n-1}).
$$
This completes the proof.
	\bigskip\bigskip
{\bf Type $\Bbb  C_n.$} There is the corresponding formula
$$
 e(\Bbb C_n) = n^2\cdot e(\Bbb A_{n-1})
$$
(with a similar proof).

	\bigskip\bigskip
{\bf Type $\Bbb D_n$.} This concerns the following diagram
$$
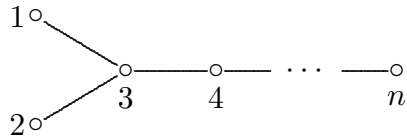
{\beginpicture
\setcoordinatesystem units <1cm,.6cm>
\multiput{$\circ$} at 0 1  0 -1  1 0  2 0  4 0 /
\put{$\cdots$} at 3 0
\plot 0.1 0.9  0.9 0.1 /
\plot 0.1 -.9  0.9 -.1 /
\plot 1.1 0  1.9 0 /
\plot 2.1 0  2.6 0 /
\plot 3.4 0  3.9 0 /
\put{$1$} at -.2 1
\put{$2$} at -.2 -1
\put{$3$} at 1 -.5
\put{$4$} at 2 -.5
\put{$n$} at 4 -.5
\endpicture}
$$
with $n \ge 4$. Actually, also the cases $n = 3$ and $n=2$ are of interest:
for $n= 3$, we have $\Bbb D_3 = \Bbb A_3$, for $n=2$ we deal with $\Bbb D_2 = \Bbb A_1\sqcup 
\Bbb A_1.$ 

Before we proceed, let us mention the following notation (see the appendix): For any $n\ge 0$, let 
$A(n) = (n+1)^{n-1}$ and $D(n) = (n\!-\!1)^n.$ 

For $k \ge 4$, we have $\Delta(k) = \Bbb D_{k-1}\sqcup \Bbb A_{n-k}$, thus the shuffling lemma yields
$$
\align
 e(\Delta(k)) &= \tbinom{n-1}{k-1} e(\Bbb D_{k-1})\cdot e(\Bbb A_{n-k}) \cr
 &= \tbinom{n-1}{k-1} 2(k-1)^k\cdot (n-k+1)^{n-k-1} \cr
 &=\tbinom{n-1}{k-1} 2D(k-1)A(n-k).
\endalign
$$ 
For $k = 3$, we have $\Delta(3) = \Bbb A_{1}\sqcup\Bbb A_1\sqcup \Bbb A_{n-3}$, and $D(2) = 1$, 
thus
$$
\align
 e(\Delta(3)) &= \tfrac{(n-1)!}{1!1!(n-3)!}e(\Bbb A_{n-3}) \cr
 &= \tbinom{n-1}2\cdot 2 \cdot (n-2)^{n-4} \cr
 &= \tbinom{n-1}2\cdot 2 D(2)A(n-3) 
\endalign
$$ 
For $k = 1$ and $k=2$, we have $\Delta(k) = \Bbb A_{n-1}$, therefore
$$
 e(\Delta(k)) = e(\Bbb A_{n-1}) = n^{n-2} = A(n-1),
$$ 
thus the sum $e(\Delta(1)) + e(\Delta(2))$ is of the form
$$
 e(\Delta(1)) + e(\Delta(2)) = \tbinom{n-1}0 2D(0)A(n-1)
$$
(since $D(0) = 1$). 

Taking into account that $D(1) = 0$, we see that
$$
\align
 \sum\nolimits_{k=1}^{n} e(\Delta(k)) &=  e(\Delta(1))+e(\Delta(2)) + \sum\nolimits_{k=3}^{n} e(\Delta(k)) \cr 
 & = \sum\nolimits_{k=1}^{n}\tbinom{n-1}{k-1} 2D(k-1)A(n-k)
\endalign
$$
but this is the coefficient $G(n\!-\!1)$ of the power series $G = D*A$, see the appendix.
The formula (3) in the appendix asserts that $G(n\!-\!1) = (n\!-\!1)^{n-1}.$

Since the Coxeter number for $\Bbb D_n$ is $h = 2(n\!-\!1),$ we have
$$
 \frac h2 \sum\nolimits_{k=1}^n e(\Delta(k)) = (n\!-\!1)\cdot 2\cdot (n\!-\!1)^{n-1} = 2(n\!-\!1)^{n},
$$
as we wanted to show.

	\bigskip\bigskip 

{\bf Type $\Bbb E_n$.} This concerns the following diagrams
$$
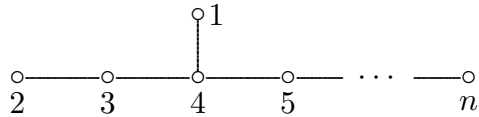
{\beginpicture
\setcoordinatesystem units <1cm,.7cm>
\multiput{$\circ$} at 0 0  1 0  2 0  3 0   2 1   5 0 /
\put{$\cdots$} at 4 0
\plot 2 0.1  2 0.9  /
\plot 0.1 0  0.9 0 /
\plot 1.1 0  1.9 0 /
\plot 2.1 0  2.9 0 /
\plot 3.1 0  3.6 0 /
\plot 4.4 0  4.9 0 /
\put{$1$} at 2.2 1
\put{$2$} at 0 -.4
\put{$3$} at 1 -.4
\put{$4$} at 2 -.4
\put{$5$} at 3 -.4

\put{$n$} at 5 -.4
\endpicture}
$$
and we will deal with the cases $n = 6,7,8$. 
	\medskip 
\noindent 
{\bf Type $\Bbb E_6$} 
$$
{\beginpicture
\setcoordinatesystem units <3cm,.6cm>
\plot -.1 -.5  2.4 -.5 /
\plot .4 0.4 .4 -4.5 /
\put{$i$} at 0 0
\put{$1$} at 0 -1
\put{$2$} at 0 -2
\put{$3$} at 0 -3
\put{$4$} at 0 -4
\put{$\Delta(i)$} at 1 0 
\put{$\Bbb A_5$} at 1 -1
\put{$\Bbb D_5$} at 1 -2 
\put{$\Bbb A_1\sqcup \Bbb A_4$} at 1 -3
\put{$\Bbb A_2\sqcup \Bbb A_1\sqcup \Bbb A_2$} at 1 -4
\put{$e(\Delta(i))$} at 2 0
\put{$1296$} at 2 -1
\put{$2048$} at 2 -2
\put{$\frac{5!}{1!4!}\; 1\cdot 125 $} at 2 -3
\put{$\frac{5!}{2!1!2!}\;3\cdot1\cdot3$} at 2 -4
\endpicture}
$$
We see:
$$
 e(\Bbb E_6) = \frac h2 \Bigl(e(\Bbb A_5)+2e(\Bbb D_5) + 2 e(\Bbb A_1\sqcup \Bbb A_4) + 
 e(\Bbb A_2\sqcup \Bbb A_1\sqcup \Bbb A_2)\Bigr) = 41\,472 = 2^93^4
$$
	\bigskip 
\noindent 
{\bf Type $\Bbb E_7$} 
$$
{\beginpicture
\setcoordinatesystem units <3cm,.6cm>
\plot -.1 -.5  2.4 -.5 /
\plot .4 0.4 .4 -7.5 /
\put{$i$} at 0 0
\put{$1$} at 0 -1
\put{$2$} at 0 -2
\put{$3$} at 0 -3
\put{$4$} at 0 -4
\put{$5$} at 0 -5
\put{$6$} at 0 -6
\put{$7$} at 0 -7
\put{$\Delta(i)$} at 1 0 
\put{$\Bbb A_6$} at 1 -1
\put{$\Bbb D_6$} at 1 -2 
\put{$\Bbb A_1\sqcup \Bbb A_5$} at 1 -3
\put{$\Bbb A_2\sqcup \Bbb A_1\sqcup \Bbb A_3$} at 1 -4
\put{$\Bbb A_4\sqcup \Bbb A_2$} at 1 -5
\put{$\Bbb D_5\sqcup \Bbb A_1$} at 1 -6
\put{$\Bbb E_6$} at 1 -7
\put{$e(\Delta(i))$} at 2 0
\put{$16807$} at 2 -1
\put{$46656$} at 2 -2
\put{$\frac{6!}{1!5!}\, 1\cdot 1296 $} at 2 -3
\put{$\frac{6!}{2!1!3!}\,3\cdot1\cdot16$} at 2 -4
\put{$\frac{6!}{4!2!}\;125\cdot3$} at 2 -5
\put{$\frac{6!}{5!1!}\,2048\cdot1$} at 2 -6
\put{$41472$} at 2 -7
\endpicture}
$$
	\medskip
$$
 e(\Bbb E_7) = 1\,062\,882= 2\cdot 3^{12}
$$
	\bigskip 
\noindent 
{\bf Type $\Bbb E_8$} 
$$
{\beginpicture
\setcoordinatesystem units <3cm,.6cm>
\plot -.1 -.5  2.4 -.5 /
\plot .4 0.4 .4 -8.5 /
\put{$i$} at 0 0
\put{$1$} at 0 -1
\put{$2$} at 0 -2
\put{$3$} at 0 -3
\put{$4$} at 0 -4
\put{$5$} at 0 -5
\put{$6$} at 0 -6
\put{$7$} at 0 -7
\put{$8$} at 0 -8
\put{$\Delta(i)$} at 1 0 
\put{$\Bbb A_7$} at 1 -1
\put{$\Bbb D_7$} at 1 -2 
\put{$\Bbb A_1\sqcup \Bbb A_6$} at 1 -3
\put{$\Bbb A_2\sqcup \Bbb A_1\sqcup \Bbb A_4$} at 1 -4
\put{$\Bbb A_4\sqcup \Bbb A_3$} at 1 -5
\put{$\Bbb D_5\sqcup \Bbb A_2$} at 1 -6
\put{$\Bbb E_6\sqcup \Bbb A_1$} at 1 -7
\put{$\Bbb E_7$} at 1 -8
\put{$e(\Delta(i))$} at 2 0
\put{$262144$} at 2 -1
\put{$559872$} at 2 -2
\put{$\frac{7!}{1!6!}\, 1\cdot 16807 $} at 2 -3
\put{$\frac{7!}{2!1!4!}\,3\cdot1\cdot125$} at 2 -4
\put{$\frac{7!}{4!3!}\;125\cdot16$} at 2 -5
\put{$\frac{7!}{5!2!}\,2048\cdot3$} at 2 -6
\put{$\frac{7!}{6!1!}\,41472\cdot 1$} at 2 -7 
\put{$1062882$} at 2 -8
\endpicture}
$$
	\medskip
$$
 e(\Bbb E_8) = 37\,968\,750 = 2\cdot 3^5\cdot 5^7.
$$

	\bigskip 
{\bf Type $\Bbb F_4$.} This concerns the following diagram
$$
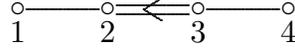
{\beginpicture
\setcoordinatesystem units <1cm,.5cm>
\multiput{$\circ$} at 0 0  1 0  2 0  3 0 /
\plot 0.1 0  0.9 0 /
\plot 1.1 0.1  1.9 0.1 /
\plot 1.1 -.1  1.9 -.1 /
\plot 2.1 0  2.9 0 /
\plot 1.6 0.25  1.4 0  1.6 -.25 /
\put{$1$} at 0 -.5
\put{$2$} at 1 -.5
\put{$3$} at 2 -.5
\put{$4$} at 3 -.5

\endpicture}
$$

$$
{\beginpicture
\setcoordinatesystem units <3cm,.6cm>
\plot -.1 -.5  2.4 -.5 /
\plot .4 0.4 .4 -4.5 /
\put{$i$} at 0 0
\put{$1$} at 0 -1
\put{$2$} at 0 -2
\put{$3$} at 0 -3
\put{$4$} at 0 -4
\put{$\Delta(i)$} at 1 0 
\put{$\Bbb B_3$} at 1 -1
\put{$\Bbb A_1\sqcup \Bbb A_2$} at 1 -2
\put{$\Bbb A_2\sqcup \Bbb A_1$} at 1 -3
\put{$\Bbb C_3$} at 1 -4 
\put{$e(\Delta(i))$} at 2 0
\put{$27$} at 2 -1
\put{$\frac{3!}{1!2!}\,1\cdot 3$} at 2 -2
\put{$\frac{3!}{2!1!}\,3\cdot 1$} at 2 -3
\put{$27$} at 2 -4
\endpicture}
$$

$$
 e(\Bbb F_4) = 432 = 2^4\cdot 3^3.
$$
	\bigskip
{\bf 6. Appendix: The binomial convolution of some power series.}
	\medskip
Let $\Bbb Z[[T]]$ be the set of formal power series $F = \sum_{n\ge 0} F(n)T^n$ 
in one variable $T$ with integer coefficients $F(n)$. 
Given power series $F = \sum_n F(n)T^n$ and $G = \sum_n G(n)T^n$,
the {\it binomial convolution} $F*G$ is by definition the power series $\sum_n H(n)T^n$
with $H(n) = \sum_k \binom n k F(k)G(n\!-\!k)$ (see [GKP]). 

We are interested in the power series $A, B, D$ with coefficients 
$A(n) = (n+1)^{n-1},$ $B(n) = n^n,$ and $D(n) = (n\!-\!1)^n,$
thus
$$
\align 
 A &= \sum_{n\ge 0} (n+1)^{n-1}T^n = 1 + T + 3T^2+16T^3 + 125 T^4+\dots \cr 
 B &= \sum_{n\ge 0} n^nT^n = 1 + T + 4 T^2 + 27 T^3 + 256 T^4 + \dots \cr
 D & = \sum_{n\ge 0} (n\!-\!1)^nT^n 
  = 1 + T^2 + 8T^3+ 81 T^4 + \dots.
\endalign
$$
The main result of the paper asserts that
$e(\Bbb A_n) = A(n)$ and $e(\Bbb B_n) = e(\Bbb C_n) = B(n)$ for $n\ge 1$ and that $e(\Bbb D_n) = 2D(n)$ 
for $n\ge 2$. 
Our proofs in section 5 
use two of the following identities, namely (1) and (3) (and we could use (2) in order to
deal with the cases $\Bbb B_n$):
	\medskip 
{\bf Proposition.}
$$
\align 
 A*A &= \sum_{n\ge 0} 2(n\!+\!2)^{n-1}T^n \tag{1}\cr 
 A*B &= \sum_{n\ge 0} (n\!+\!1)^{n}T^n \tag{2}\cr 
 A*D &= \sum_{n\ge 0} n^{n}T^n  = B \tag{3}
\endalign
$$
	\medskip
Proof. Let us recall Abel's identity [Ab]
$$
 (x+y)^n = \sum_{k=0}^n \binom n k x(x-kz)^{k-1}(y+kz)^{n-k}
$$
which is valid in any commutative ring with $x$ being invertible.
Several proofs can be found in Comptet [Co]. We need Abel's identity for $x=1$ and $z=-1$, thus
the identity
$$
 (1+y)^n = \sum_{k=0}^n \binom n k (1+k)^{k-1}(y-k)^{n-k}.
$$

Let us start with the proof of (2), using Abel's identity
for  $y = n$ (and $x=1,\ z = -1$):
$$
  (1\!+\!n)^{n} =  \sum_{k=0}^n \binom n k (1+k)^{k-1}(n-k)^{n-k} 
  = \sum_{k=0}^n \binom n k A(k)B(n-k) = (A*B)(n).
$$
For the proof of (3), we use Abel's identity for $y = n-1$ (and $x=1,\ z = -1$):
$$
\align
  n^{n} &= (1+(n-1))^{n} 
  =  \sum_{k=0}^n \binom n k (1+k)^{k-1}(n-1-k)^{n-k} \cr
  &= \sum_{k=0}^n \binom n k A(k)D(n-k) = (A*D)(n).
\endalign
$$
For the proof of (1), we expand $(n+2)^{n-1}$ with $y = n+1$ (and again $x=1,\ z = -1$):
$$
  (1+(n+1))^{n-1} 
   = \sum_{k=0}^{n}\binom{n-1}k (1+k)^{k-1}(n+1-k)^{n-1-k},  \tag{$*$}
$$
note that we have added the summand with index $k=n$; there is no harm, since 
by definition $\binom n{n-1} = 0.$ Replacing the summation index $k$ by $n-k$, 
and using the equality $\binom{n-1}{n-k} = \binom{n-1}{k-1}$, we see that we also have
$$
 (1+(n+1))^{n-1} = \sum_{k=0}^{n}\binom{n-1}{k-1} (1+n-k)^{n-k-1}(k+1)^{k-1}. \tag{$**$}
$$
Since $\binom{n-1}k+\binom{n-1}{k-1} = \binom nk,$ the summation of $(*)$ and $(**)$ yields
$$
\align
 2(n+2)^{n-1} &= \sum_{k=0}^{n}\binom{n}{k} (k+1)^{k-1}(n-k+1)^{n-k-1} \cr
    &= \sum_{k=0}^{n}\binom{n}{k} A(k)A(n-k) = (A*A)(n).
\endalign
$$
This completes the proof of the Proposition.
	\medskip
It seems to us that these binomial convolution formulas are very pretty;
as an example, let us exhibit the coefficients of $T^4$ in $A*A,\ A*B,\ A*D$:
$$
{\beginpicture
\setcoordinatesystem units <.9cm,.5cm>
\multiput{$=$} at 11.2 0   11.2 -1  11.2 -2 /

\multiput{$+$\strut} at 3 0  5 0  7 0  9 0 /
\multiput{$+$\strut} at 3 -1  5 -1  7 -1  9 -1 /
\multiput{$+$\strut} at 3 -2  5 -2  7 -2  9 -2 /
\put{$2\cdot 6^3$\strut} [r] at 12.5 0
\put{$1 \cdot 1\cdot 125$\strut} at 2 0
\put{$4 \cdot 1\cdot 16$\strut} at 4 0
\put{$6 \cdot 3\cdot 3$\strut} at 6 0
\put{$4 \cdot 16\cdot 1$\strut} at 8 0
\put{$1 \cdot 125\cdot 1$\strut} at 10 0
\put{$5^4$\strut} [r] at 12.5 -1
\put{$1 \cdot 1\cdot 256$\strut} at 2 -1
\put{$4 \cdot 1\cdot 27$\strut} at 4 -1
\put{$6 \cdot 3\cdot 4$\strut} at 6 -1
\put{$4 \cdot 16\cdot 1$\strut} at 8 -1
\put{$1 \cdot 125\cdot 1$\strut} at 10 -1
\put{$4^4$\strut} [r] at 12.5 -2
\put{$1 \cdot 1\cdot\ \, 81$\strut} at 2 -2
\put{$4 \cdot 1\cdot\  8$\strut} at 4 -2
\put{$6 \cdot 3\cdot 1$\strut} at 6 -2
\put{$4 \cdot 16\cdot 0$\strut} at 8 -2
\put{$1 \cdot 125\cdot 1$\strut} at 10 -2
\put{$(A*A)(4)$\strut} at -.5 0 
\put{$(A*B)(4)$\strut} at -.5 -1 
\put{$(A*D)(4)$\strut} at -.5 -2
\endpicture}
$$
	\bigskip 
Finally, let us add some general information
concerning the sequences $A,B,D$ as provided by  
Sloane's On-Line Encyclopedia of Integer Sequences [Sl].
The sequence $A(n) = (n+1)^{n-1}$ is the Sloane sequence A000272, 
but shifted by $1$, thus $A(n)$
is the {\it number of trees on $n+1$ labeled nodes.}
The sequence $B(n) = n^n$ is the Sloane sequence A000312, the number $B(n)$
is the 	
{\it number of functions from the set $\{1,2,...,n\}$ to itself.}
The sequence $D(n) = (n\!-\!1)^n$ with $e(\Bbb D_n) = 2D(n)$ for $n\ge 2$ is the Sloane sequence A065440;
the number $D(n)$ is the {\it number of functions from the set 
$\{1,2,...,n\}$ to itself 
without fixed points.}

Here are the first terms of the sequences $A, B, 2D$, namely 
$A(n), B(n), 2D(n),$ with $n \le 10$; note that
$A(n) = e(\Bbb A_n), B(n) = e(\Bbb B_n)$, for $n\ge 1$ and $2D(n) = e(\Bbb D_n)$, for $n\ge 2$.

$$
{\beginpicture
\setcoordinatesystem units <3cm,.5cm>
\put{$n$} at -.3 1
\put{$A(n)$} at .7 1
\put{$B(n)$} at 1.7 1
\put{$2D(n)$} at 2.7 1
\put{$0$} at -.3 0
\put{$1$} at -.3 -1
\put{$2$} at -.3 -2
\put{$3$} at -.3 -3
\put{$4$} at -.3 -4
\put{$5$} at -.3 -5
\put{$6$} at -.3 -6
\put{$7$} at -.3 -7
\put{$8$} at -.3 -8
\put{$9$} at -.3 -9
\put{$10$} at -.3 -10
\put{$1$} [r] at 1 0
\put{$1$} [r] at 1 -1
\put{$3$} [r] at 1 -2
\put{$16$} [r] at 1 -3
\put{$125$} [r] at 1 -4
\put{$1\,296$} [r] at 1 -5
\put{$16\,807$} [r] at 1 -6
\put{$262\,144$} [r] at 1 -7
\put{$4\,782\,969$} [r] at 1 -8
\put{$100\,000\,000$} [r] at 1 -9
\put{$2\,357\,947\,691$} [r] at 1 -10

\put{$1$} [r] at  2 0
\put{$1$} [r] at  2 -1
\put{$4$} [r] at  2 -2
\put{$27$} [r] at  2 -3
\put{$256$} [r] at  2 -4
\put{$3\,125$} [r] at  2 -5
\put{$46\,656$} [r] at  2 -6
\put{$823\,543$} [r] at  2 -7
\put{$12\,777\,216$} [r] at  2 -8
\put{$387\,420\,489$} [r] at  2 -9
\put{$10\,000\,000\,000$} [r] at  2 -10
\put{$2$} [r] at  3 0
\put{$0$} [r] at  3 -1
\put{$2$} [r] at  3 -2
\put{$16$} [r] at  3 -3
\put{$162$} [r] at  3 -4
\put{$2\,048$} [r] at  3 -5
\put{$31\,250$} [r] at  3 -6
\put{$559\,872$} [r] at  3 -7
\put{$11\,529\,602$} [r] at  3 -8
\put{$268\,435\,456$} [r] at  3 -9
\put{$6\,973\,568\,802$} [r] at  3 -10

\endpicture}
$$

	\bigskip

{\bf 7\. References.}
		     \medskip
\item{[Ab]} N\. Abel: Beweis eines Ausdrucks, von welchem die Binomial-Formel ein
	einzelner Fall ist. Crelle's J\. Math. 1 (1826), 159-160.
\item{[Ag]} D\. Armstrong: Generalized Noncrossing Partitions and
	Combinatorics of Coxeter Groups. 
        Memoirs of the Amer\. Math\. Soc\. 949 (2009).
\item{[AR]} C\. A\. Athanasiadis, V\. Reiner: Noncrossing partitions
	 for the group $D_n.$ SIAM J. Discrete Math. 18 (2004), no. 2, 397-417.
\item{[B]} N\. Bourbaki: Groupes et algebres de Lie: 
    Chapitres 4, 5 et 6. Paris (1968). 
\item{[Ch]} F\. Chapoton: Enumerative properties of generalized associahedra.
  Sem\. Lothar\. Combin. 51 (2004/5).
\item{[CS]} G\. Chapuy, C\. Stump: Counting factorizations of Coxeter elements
   into products of reflections. arXiv:1211.2789.
\item{[Co]} L\. Comptet: Advanced Combinatorics. Reidel (1974).
\item{[CB]} W\. Crawley-Boevey: Exceptional sequences of quivers. In:
    Canadian Math\. Soc\. Proceedings 14 (1993), 117-124.
\item{[D]} P\. Deligne: Letter to E\. Looijenga 9.3.1974. Online available:\newline
   http://homepage.univie.ac.at/christian.stump/Deligne\_Looijenga\_Letter\_09-03-1974.pdf
\item{[DR1]} V\. Dlab, C\. M\. Ringel: On algebras of finite representation type. J\. Algebra 33
      (1975), 306-394.
\item{[DR2]} V\. Dlab, C\. M\. Ringel: Indecomposable representations of graphs and algebras.
      Mem\. Amer\. Math\. Soc\. 173 (1976).
\item{[DRS]} P\. Dowbor, C\. M\. Ringel, D\. Simson:
        Hereditary artinian rings of finite representation type. Proceedings ICRA 2. Springer LNM 832
       (1980),    232-241.
\item{[FR]} W\. Fakieh, C\. M\. Ringel: The hereditary artinian rings of type $\Bbb H_3$
   and   $\Bbb H_4$. In preparation. 
\item{[GL]} W\. Geigle, H\. Lenzing: Perpendicular categories with applications to representations 
     and  sheaves, J. Algebra 144 (1991), 273-343
\item{[GKP]} R\. L\. Graham, D\. E\. Knuth, O\. Patashnik: Concrete Mathematics. A Foundation for
        Computer Science. Addison Wesley, Reading (1989).
\item{[HK]} A\. Hubery,
      H\. Krause: A categorification of noncrossing partitions. In preparation.
\item{[IT]} C\. Ingalls, H\. Thomas: Noncrossing partitions and representations of quivers.
      Comp\. Math\. 145 (2009), 1533-1562.
\item{[K]} G\. Kreweras: Sur les partitions non crois\'ees d'un cycle.
      Discr. Math. 1, number 4 (1972), 333-350 
\item{[L]}  E\. Looijenga: The complement of the bifurcation
      variety of a simple singularity. Invent. Math. 23 (1974), 105-116.
\item{[O]}S\. Oppermann: Auslander-Reiten theory of representation directed artinian rings.
      Diplomarbeit, Stuttgart 2005. http://www.math.ntnu.no/$\sim$opperman/artinian.pdf
\item{[Rd]} N\. Reading: Chains in the noncrossing partition lattice.
      SIAM J. Discrete Math. 22 (2008), no. 3, 875-886.
\item{[Rn]} V\. Reiner, Non-crossing partitions for classical reflection
      groups. Discrete Math. 177 (1997), no. 1-3, 195-222.
\item{[R1]} C\. M\. Ringel: Tame algebras and integral quadratic forms. Springer LNM 1099 (1984).
\item{[R2]} C\. M\. Ringel: The braid group action on the set of exceptional sequences of a
      hereditary    algebra.
          In: Abelian Group Theory and Related Topics. Contemp\. Math\. 171 (1994), 339-352.
\item{[R3]} C\. M\. Ringel: Some remarks concerning tilting modules and tilted algebras.
      Origin. Relevance. Future. (An appendix to the Handbook of Tilting Theory.)
      London Math\. Soc\. Lecture Note Series 332. Cambridge University Press (2007),
      413-472.
\item{[Se]} U\. Seidel: Exceptional sequences for quivers of Dynkin type. Comm\. Algebra
      29 (2001). 1373-1386. 
\item{[S1]} A\. Schofield: 
      Hereditary artinian rings of finite representation type
      and extensions of simple artinian rings. Math. Proc. Cambridge Philos. Soc.
      102 (1987), no. 3, 411-420.

\item{[S2]} A\. Schofield: Semi-invariants of quivers. J. London Math. Soc. (2) 43 (1991), 385–395.
\item{[Sl]} N\. J\. A\. Sloane: On-Line Encyclopedia of Integer Sequences. http://oeis.org/
	\bigskip\bigskip 

{\rmk
\noindent 
Mustafa A. A. Obaid: E-mail: {\ttk drmobaid\@yahoo.com}
	\smallskip
\noindent 
S. Khalid Nauman: E-mail: {\ttk snauman\@kau.edu.sa}
	\smallskip
\noindent 
Wafa S. Al Shammakh: E-mail: {\ttk wafasalem25\@hotmail.com}
	\smallskip
\noindent 
Wafaa  M. Fakieh: E-mail: {\ttk wafaa.fakieh\@hotmail.com}
	\smallskip
\noindent 
Claus Michael Ringel: E-mail: {\ttk ringel\@math.uni-bielefeld.de}
	\medskip
\noindent 
King Abdulaziz University, Faculty of Science, \par 
P.O.Box 80203, Jeddah 21589, Saudi Arabia

}

\bye